\documentclass[11pt]{article}

\usepackage{a4}
\usepackage{amsmath}
\usepackage{amssymb}
\usepackage[PostScript=dvips]{diagrams}

\DeclareMathAlphabet{\scr}{U}{eus}{m}{n}

\newcommand\N{{\mathbb N}}
\newcommand\Z{{\mathbb Z}}
\newcommand\Q{{\mathbb Q}}

\newcommand\cF{{\cal F}}

\newcommand\cO{{\cal O}}

\newcommand\Mn{\mathbf M^{\rm naive}}
\newcommand\M{\mathbf M}

\newcommand\Mloc{{\mathbf M}^{\rm loc}}
\newcommand\Mcan{{\mathbf M}^{\rm can}}

\newcommand\tens\otimes

\newcommand\lto{\longrightarrow}

\renewcommand{\hom}{\mathop{\rm Hom}\nolimits}

\newcommand{\Spec}{\mathop{\rm Spec}}
\newcommand{\Gal}{\mathop{\rm Gal}}
\newcommand{\diag}{\mathop{\rm diag}}

\newcommand{\Res}{\mathop{\rm Res}\nolimits}
\renewcommand{\det}{\mathop{\rm det}\nolimits}

\newcommand{\Tr}{\mathop{\rm Tr}}
\newcommand{\Perm}{\mathop{\rm Perm}\nolimits}
\newcommand{\Adm}{\mathop{\rm Adm}\nolimits}
\newcommand{\faces}{\mathcal F}
\newcommand{\Waff}{W_{\rm aff}}

\newcommand\qed{\hfill$\square$}

\newtheorem{thm}{Theorem}[section]
\newtheorem{stz}[thm]{Proposition}
\newtheorem{lem}[thm]{Lemma}
\newtheorem{Def}[thm]{Definition}
\newtheorem{kor}[thm]{Corollary}

\newtheorem{conj}[thm]{Conjecture}

\begin{document}

\title{Topological flatness of local models \\in the ramified case}
\author{Ulrich G\"ortz \footnote{Mathematisches Institut, Universit\"at zu K\"oln,
Weyertal 86--90, 50931 K\"oln, Germany. Email: ugoertz@mi.uni-koeln.de}}
\date{}
\maketitle

{\bf Abstract.}
Local models are schemes defined in terms of linear algebra which can 
be used to study the local structure of integral models of certain Shimura
varieties, with parahoric level structure. We investigate the local models
for groups of the form $\Res_{F/\Q_p} GL_n$ and $\Res_{F/\Q_p} GSp_{2g}$ where
$F/Q_p$ is a totally ramified extension, as defined by Pappas and Rapoport,
and show that they are topologically flat. In the linear case, flatness can
be deduced from this. 

\vskip.4cm

\section{Introduction}
\markboth{U. G\"ortz, Topological flatness of local models}
{U. G\"ortz, Topological flatness of local models}

It is an interesting problem to define 'good' models of Shimura 
varieties over the ring of integers of the reflex field, or
over its completion at some prime ideal. For Shimura varieties of
PEL type, which can be described as moduli spaces of abelian varieties,
it is desirable to define such a model by a moduli problem, too.

In their book \cite{RZ}, Rapoport and Zink define such models 
in the case of parahoric level structures. They also define
the so-called local models which \'etale-locally around each point
of the special fibre are isomorphic to the model of the Shimura
variety, but which are defined purely in terms of linear algebra.
They provide a very useful tool to investigate local properties
of the corresponding models.

Rapoport and Zink conjectured that these models are flat, which 
is a property a
good model certainly should have. The conjecture is true for 
Shimura varieties associated to unitary or symplectic groups
that split over an unramified extension of $\Q_p$ (see \cite{G},
\cite{G2}). But Pappas \cite{P} showed that it is in general 
false if the group splits only after a ramified base change.
Often the models are not even topologically flat, i. e.
the closure of the generic fibre is not even set-theoretically
equal to the whole model. We follow the terminology of Pappas and 
Rapoport and call the local model associated to these models of the
Shimura variety the {\em naive local model}.

For groups of the form $\Res_{F/\Q_p}GL_n$, where $F/\Q_p$ is a 
totally ramified extension,
Pappas and Rapoport \cite{PR} then defined a new local model, which they
call the {\em local model}. In the case 
of level structures corresponding to a maximal parahoric subgroup,
it is essentially defined as the closure of the
generic fibre in the old local model, so it is automatically flat. 
Pappas and Rapoport showed that
this new model has several good properties; for 
example its special fibre is normal, with Cohen-Macaulay singularities, 
and over a finite extension of $\cO_E$ admits a resolution of singularities.
The definition of the local model in the maximal parahoric case
gives rise to a definition of a new local model in the general
parahoric case. It is not obvious anymore that these general local 
models are flat, too, and it is the purpose of this note to show that
they are topologically flat. Together with the methods
used in \cite{G} (Frobenius splitting of Schubert varieties in the affine 
flag variety), one can then infer flatness. We therefore obtain 
the following theorem:

\begin{thm}
Let $F/\Q_p$ be a (possibly ramified) finite extension. Then
the local model associated to $\Res_{F/\Q_p}GL_n$ is flat
over $\Z_p$.
\end{thm}

An essential ingredient of the proof is a combinatorial result of Haines
and Ng\^o about the so-called $\mu$-admissible and $\mu$-permissible sets
(see below for further details). This result, as it stands, relates to the
Iwahori case. In section \ref{parahoric} we show how it can be used to
prove the corresponding result in the general parahoric case.

In section \ref{symplectic} we show that the local
model for groups of the form $\Res_{F/\Q_p}GSp_{2g}$,
 $F/\Q_p$ totally ramified, is topologically flat, too.
(In this case the naive local model and the local
model coincide topologically.) But here it is more difficult to
deduce the flatness from the topological flatness, and we can only
make a conjecture in this direction. On the other hand, it seems
reasonable to expect that in this case the naive local model is
flat itself (and thus coincides with the local model).

In a subsequent article \cite{PR2} Pappas and Rapoport define 
yet another model, which they call the {\em canonical model}.
Loosely speaking, it is the image of a certain morphism 
from a twisted product of unramified local models to the naive local 
model. In particular it is always flat, because the unramified
local models are flat, and in addition has other good properties.
Nevertheless the 
flatness question for the original local model remains interesting; 
if the local model is flat, it must coincide with the canonical
local model. On the one hand this shows that the special fibre of the 
canonical local model has the 'correct' combinatorial description,
on the other hand one gets an interesting 'resolution' of the local
model.

The results of Pappas and Rapoport together with the
combinatorial results of Haines and Ng\^o, and in section \ref{parahoric},
respectively, can be used to give a shorter (but less elementary) proof of
the topological flatness, as came out of discussions with Haines.
See the remark after proposition \ref{prop_flatness} for an outline.

I would like to thank T. Haines for many fruitful conversations on
this topic.
I am also grateful to G. Pappas, M. Rapoport and T. Wedhorn for their
helpful remarks.
The main part of this article was written during my stay at the 
Institute for Advanced Study in Princeton. I would like to thank the
Institute for providing a great working environment. The support of the
Deutsche Forschungsgemeinschaft and the National Science Foundation (grant
no. DMS 97-29992) is gratefully acknowledged.

\section{Definitions}

We will use the notation of Pappas and Rapoport \cite{PR}. 
Let us repeat part of it:

Let $F_0$ be a field, complete with respect to a non-archimedean 
valuation. Let $\cO_{F_0}$ denote its ring of integers, $\pi_0$
a uniformizer. We assume that the residue class field is perfect.
We fix a separable closure $F_0^{\rm sep}$ of $F_0$.

Let $F/F_0$ be a totally ramified extension of degree $e$ inside $F_0^{\rm sep}$.
Let $\cO_F$ be the ring of integers of $F$, and let $\pi\in\cO_F$ be
a uniformizer which is the root of an Eisenstein polynomial
$T^e + \sum a_k T^k$. For each embedding $F\lto F_0^{\rm sep}$
we choose an integer $0\le r_\varphi \le d$.
Associated to these data we have the reflex field $E$, a finite 
extension of $F_0$ contained in $F_0^{\rm sep}$, which is defined
by $\Gal(F_0^{\rm sep}/E) = \{ \sigma \in \Gal(F_0^{\rm sep}/F_0);\ 
 \forall \varphi: \ r_{\sigma\varphi} = r_\varphi \}$.

Further, let $V = F^d$, $\Lambda_0 = \cO_F^d \subset V$, and 
denote the canonical basis by $e_1, \dots, e_d$.

Choose $I \subseteq \{1, \dots, d \}$.

Let us recall the definition of the 'naive' local model
 $\Mn = \M(e,d,(\Lambda_i)_{i\in I}, (r_\varphi)_\varphi)$ 
(of course $\Mn$ depends
on $F/F_0$, not just on $e$, so this is an abuse of notation
which nevertheless seems useful):

It is defined over the ring of integers $\cO_E$ of the
reflex field $E$, and
its $S$-valued points are the isomorphism classes
of commutative diagrams 
\begin{diagram}
\Lambda_{i_0, S} & \rTo & \Lambda_{i_1,S} & \rTo & \cdots 
& \rTo & \Lambda_{i_{m-1},S} & \rTo^\pi & \Lambda_{i_0,S} \\
\uInto              &      & \uInto        &      &   &&\uInto & &\uInto   \\
\cF_0     & \rTo & \cF_1  & \rTo & \cdots & \rTo & \cF_{m-1}
  &\rTo & \cF_0
\end{diagram}
where $\Lambda_{i,S}$ is $\Lambda_i \tens_{O} \cO_S$, and where 
the $\cF_\kappa$ are 
 $\cO_F \tens_{\cO_{F_0}} \cO_S$-submodules. We require that locally
on $S$, the $\cF_i$ are direct summands as $\cO_S$-modules and
that we have the following identity of polynomials
('determinant condition'):
$$ \det_{\cO_S}(T - \Pi|\cF_i) = \prod_\varphi (T-\varphi(\pi))^{r_\varphi}. $$

Pappas and Rapoport show that this naive local model is almost never flat,
even if $I$ consists of only one element. In this case, i. e. when $I$
consists of only one element, they define
a new local model $\Mloc$ as the scheme-theoretic closure of the generic fibre
in the naive local model, and show that this new local model 
has several good properties. In particular, its special fibre is reduced, 
and can be described as a union of Schubert varieties with the 'right'
index set.

Based on this, they define a new Iwahori type local model
$\Mloc$ as the closed
subscheme of the naive local model such that all projections to
the parahoric local model map to the new local model.
Below we will show that this local model is topologically flat.
Furthermore, using the technique of Frobenius splittings, one can
see that its special fibre is reduced. Thus $\Mloc$ is flat.

\section{A certain map between local models}
\label{certainmap}

It is clear that in the definition of the local model,
we have not used the fact that the sequence $\Lambda_0 \lto \cdots
 \lto \Lambda_{n-1} \lto \pi^{-1}\Lambda_0$ is a lattice chain --- 
in fact we could
make a completely analogous definition for any sequence of free $\cO_F$-modules.
(Compare \cite{G}, 'general schemes of compatible subspaces'.)

In particular, consider the following situation:
Choose a partition $\{0, \dots, d-1 \} = \coprod_\alpha I_\alpha$.
We can then decompose the lattices $\Lambda_i$ as
$$ \Lambda_i = \bigoplus_\alpha \Lambda_i^\alpha, $$
where $\Lambda_i^\alpha = \bigoplus_{j\in I_\alpha} \cO_F e^i_j$.
Here the $e^i_j = \left\{ \begin{array}{ll} \pi^{-1}e_j & j \le i \\
 e_j & j > i \end{array}\right.$ denote the canonical $\cO_F$-basis 
of $\Lambda_i$. 

Now choose $0 \le r_\varphi^\alpha \le |I_\alpha|$ such that 
$\sum_\alpha r_\varphi^\alpha = r_\varphi$.

This gives rise to a decomposition of our lattice chain. More precisely
for each $\alpha$ we get a sequence
$$ \Lambda_0^\alpha  \lto \cdots
 \lto \Lambda_{n-1}^\alpha \lto \pi^{-1}\Lambda_0^\alpha $$
of free $\cO_F$-modules. The difference  between these sequences and the
original lattice chain (apart from the rank of
the lattices) is that now some of the maps may be the identity.
Nevertheless, we can define a 'local model' 
 $\M^\alpha = \M(e,d,(\Lambda_i^\alpha)_{i=0, \dots, n-1}, 
                   (r_\varphi^\alpha)_\varphi)$. 
It is isomorphic to some local model in the original sense (because as far as
the local model is concerned, we can just omit the 'identity steps').

The reflex field associated to a tuple $(r_\varphi^\alpha)_\varphi$
will in general be different from the one belonging to 
 $(r_\varphi)_\varphi$. Thus the $\M^\alpha$ will in general not 
be defined over the same ring as the naive local model $\M/\cO_E$ associated
to the $r_\varphi$. To simplify the situation, we make the following
additional assumption: whenever $r_\varphi = r_\psi$, we have
 $r_\varphi^\alpha = r_\psi^\alpha$ for all $\alpha$. Under this 
assumption, all $\M^\alpha$ will be defined over $\cO_E$, or even over
a smaller ring. If necessary we apply a base change, and in the 
following we consider all the $\M^\alpha$ as $\cO_E$-schemes.

We then have a canonical map
\begin{eqnarray*}
\prod_\alpha \M(e,d,(\Lambda_i^\alpha)_{i=0, \dots, n-1}, 
   (r_\varphi^\alpha)_\varphi)
& \lto & \M(e,d,(\Lambda_i)_{i=0, \dots, n-1}, (r_\varphi)_\varphi) \\
((\cF^\alpha_i)_i)_\alpha) & \mapsto & \left(\bigoplus_\alpha \cF_i^\alpha\right)_i.
\end{eqnarray*}

\section{The case $d=1$}

In this section, we will look at the (trivial) case where $d=1$.
We will use the following observations in section \ref{lifting}.

Let us give a description of the local model corresponding
to a maximal parahoric subgroup (the standard local model
in the sense of \cite{PR}, section 2). (Actually, in this case
this coincides with the Iwahori type local model.)

So, we have integers $0 \le r_\varphi\le 1$, $\varphi = 1, \dots, e$.
Let us assume, for notational convenience, that the $r_\varphi$ are
in descending order, such that they are completely determined by
 $r = \sum r_\varphi$. 

The local model in this case is just
$\Spec \cO_E$, where $E$ is the corresponding reflex field. 
The point in the special fibre is the subspace $\cF$
corresponding to the $e \times r$-matrix
$$ \overline{M} = \left( \begin{array}{c} I_r \\ 0 \end{array}
\right);
$$
the operator $\Pi|\cF$ has Jordan type $(r)$.
Let us denote the subspace corresponding to the $R$-valued point
of the local model
by $\cF(e,r)$.
The description of this subspace in
terms of matrices is the following. We have a matrix
$$ M(e,r) = M = \left( \begin{array}{c} I_r \\ 
    (b_{ij})_{{ i=1, \dots, e-r \atop 
                               j = 1, \dots r }}
\end{array}
\right),
$$
such that 
$$ \Pi M = M A,$$
for some matrix $A$ with characteristic polynomial
$$ \det(T-A) = \prod (T-\varphi(\Pi))^{r_\varphi}. $$
Furthermore, the reduction of $M$ modulo $\pi$ is $\overline{M}$.

{\bf Example.} Assume that $\pi = \pi_0^6$. Now, if we let $r=2$
and choose the $r_\varphi$ such that the resulting characteristic 
polynomial is $T^2 - \pi^2$, we get
$$ M(6,2) = \left( \begin{array}{cc}
1 & \\
 & 1 \\
\pi^2 & \\
 & \pi^2 \\
\pi^4 & \\
 & \pi^4 
\end{array} \right). $$

\section{Lifting of points}
\label{lifting}

\begin{stz} \label{prop_flatness}
The local model $\Mloc$ is topologically flat, i. e.
the generic points of the irreducible
components can be lifted to the special fibre.
\end{stz}

{\em Remark.}
The following easier, though less elementary, proof resulted 
from discussions with Haines.
It uses the theory of the splitting model developed by
Pappas and Rapoport.

In \cite{PR2}, Pappas and Rapoport define the {\em splitting model},
which is a twisted product of unramified local models, and which maps to
the local model $\Mloc$.
By definition, the {\em canonical local model} $\Mcan$ is the image of 
this morphism. 

It is not hard to see that the special fibre of $\Mcan$
consists of the Schubert cells corresponding to the elements of the
$\mu$-admissible set. By the theorem of Haines and Ng\^o (and the
generalization in section \ref{parahoric}, respectively), this set 
coincides with the $\mu$-permissible set, which parametrizes the Schubert
cells in the special fibre of $\Mloc$.

Since the splitting model is flat by \cite{G}, the
canonical local model is flat, and thus the local model $\Mloc$ is
topologically flat.

We now give a more direct proof of the proposition.

{\em Proof.}
In this section, we consider only the Iwahori case. It will be obvious
though, that the proof carries over to the general parahoric case without
any problems once the results of Haines and Ng\^o are generalized
correspondingly. In section \ref{parahoric} we will give a proof of the
more general statement which we need.

We can embed the special fibre of $\Mloc$ into the affine flag variety
for $GL_d$ in the standard way, see \cite{G}.

By the definition of the local model and since the
special fibers of the maximal-parahoric new local models have 
the right stratification,
the special fibre is the union of the Schubert cells corresponding to
the $\mu$-permissible alcoves.
Here $\mu = (\mu_1, \dots, \mu_d)$ is the dual partition to 
 $(r_\varphi)_\varphi$.

In order to show that the local model associated to some dominant coweight
 $\mu$ is topologically flat, we need to
have a good understanding of the set of irreducible components of its
special fibre. This amounts to a combinatorial problem in the extended
affine Weyl group, namely to relating the so-called $\mu$-permissible and
$\mu$-admissible sets.

By recent work of Haines and Ng\^o \cite{HN} we know that in this case
(as in the minuscule case), the set of $\mu$-permissible alcoves
coincides with the set of $\mu$-admissible alcoves. Thus the
maximal elements (with respect to the Bruhat order), which correspond
to the irreducible components of the special fibre, are just 
the conjugates of $\mu$ under the finite Weyl group, i.e. the
permutations of $\mu$ under $S_d$.

Choose a permutation $\nu$ of $\mu$ and
write $\nu = (\nu_1, \dots, \nu_d)$.

A key observation is the following lemma which came up in a 
discussion with T. Haines.

\begin{lem} \label{keylemma}
Let $R$ be a discrete valuation ring, let $s$ be the closed
and $\eta$ the generic point of $\Spec R$. Let $U$ be an $R$-scheme
of finite type,
and $x \in U_s$ a closed point of the special fibre which lies
in a unique irreducible component $S$ of $U_s$.

Furthermore, assume that $U_\eta$
is irreducible, that $\dim S_s = \dim U_\eta$, and that
 $x$ can be lifted to a closed point of $U_\eta$.

Then the generic point of $S$ can be lifted to $U_\eta$.
\end{lem}

{\em Proof.} Denote by $U'$ the scheme-theoretic closure
of $U_\eta$ in $U$. Then $U'$ is flat over $R$, and $x \in U'_s$.
Furthermore, $\dim \cO_{U',x} = \dim \cO_{U'_s,x} + 1$ (see
\cite{M}, theorem 15.1), and $\dim \cO_{U',x} \ge \dim U_\eta +1$,
since $x$ can be lifted to a closed point of $U_\eta$.

Thus $\dim \cO_{U'_s,x} \ge \dim U_\eta = \dim S$, which implies
that $S \subseteq U'$. \qed

Although the lemma is easy to prove, it reduces the combinatorial
difficulty of our task enormously.

The lemma shows that it really is enough to show
that one suitably chosen closed point of the special fibre
can be lifted to the generic fibre. We will choose the point 
 $x$ = $(\cF_i)_i$, where each $\cF_i$ is given
by the $de \times r$-matrix

$$\overline{M} =  \left( \begin{array}{cccc}  
\begin{array}{c} I_{\nu_1} \\
0_{e-\nu_1 \times \nu_1}\end{array} & &  &\\
  &   \begin{array}{c} I_{\nu_2} \\
0_{e-\nu_2 \times \nu_2}\end{array} & & \\
  &  & \ddots & \\
  &  & & \begin{array}{c} I_{\nu_d} \\
0_{e-\nu_d \times \nu_d}\end{array}
\end{array} \right). $$

It is clear that $x$ lies in the stratum corresponding to $\nu$. 
(Look at the Jordan type of $\Pi|\cF_i$.)

Now how can we lift this point to the
generic fibre?
We will describe two different methods to do this. The first one is more
elementary, the second one more conceptual. 

\subsection{First method (elementary)}

To simplify the notation, let us assume that the $r_\varphi$ are in 
descending order.
 We will build the corresponding matrix (over $\cO_E$) 
by putting together several matrices of the form $M(e,s)$.
More precisely, consider the $\cO_E$-valued point where each subspace
 $\cF_i$ is given by the same matrix 

$$M =  \left( \begin{array}{cccc} 
 M(e,\nu_1) & & &\\
  &    M(e,\nu_2) & & \\
  & &  \ddots & \\
  & &  & M(e,\nu_d)
\end{array} \right). $$

Note that the coefficients of this matrix are indeed contained in 
 $\cO_E$.
Clearly, the reduction of $M$ is $\overline{M}$. We have to check that

\begin{itemize}
\item $M$ describes a subspace that lies in the maximal-parahoric
local model (determinant condition)

\item The subspace $\cF_i$ is mapped to $\cF_{i+1}$.
\end{itemize}

But both conditions are clearly satisfied, since the $M(e,\nu_i)$
satisfy the corresponding conditions and since there are no
interactions between the blocks. Let us make that more precise.

To check the first condition, fix any $i$. Note that $\cF_i$
is $\Pi$-invariant: this means that there is a matrix $A$
such that 
$$ \diag(\Pi, \dots, \Pi)M = MA, $$
and we can simply take $A = \Pi|\cF = \diag(\Pi|\cF(e,\nu_1),\dots,
\Pi|\cF(e,\nu_d))$.
Now let us check the determinant condition;
the characteristic polynomial of $\Pi|\cF$ is just the
product of the characteristic polynomials of the $\Pi|\cF(e,\nu_i)$,
so we get 
$$ 
\prod_{i=1}^d \prod_{\varphi = 1}^{\nu_i} T-\varphi(\Pi)
= \prod_{\varphi=1}^e (T-\varphi(\Pi))^{r_\varphi}, $$
as it should be.

The second condition is easily checked, too. We just have to observe
that 
$$ \diag(1, \dots, 1, \Pi, 1, \dots, 1) M = M \diag(1, \dots,
   1, \Pi|\cF(e,\nu_i), 1, \dots, 1).
$$

Thus we have indeed found a lifting of our point to the generic fibre.

\subsection{Second method (conceptual)}

In this section, we will describe the $R$-valued point constructed in the
previous section as a morphism from a product of trivial (i.e. $\cong \Spec \cO_E$)
local models to the given local model, which comes from the situation studied
in section \ref{certainmap}.

As partition of $I$, we will choose the partition into singleton sets:
 $I_\alpha = \{ \alpha \}$, $\alpha = 0, \dots, d-1$.
Furthermore we choose $0\le r_\varphi^\alpha \le 1$ such that
$$ \sum_\alpha r_\varphi^\alpha = r_\varphi, \quad \text{and } 
   \sum_\varphi r_\varphi^\alpha = \nu_\alpha. $$
It is easy to see that such $r_\varphi^\alpha$ exist and that they are
uniquely determined. Furthermore, whenever $r_\varphi = r_\psi$, we have
 $r_\varphi^\alpha = r_\psi^\alpha$ for all $\alpha$. Thus all
 $\M(e,d,(\Lambda_i^\alpha)_{i=0, \dots, n-1}, (r_\varphi^\alpha)_\varphi)$
are defined over a subring of $\cO_E$, and we denote the base change
to $\cO_E$ by $\M^\alpha$.

Now consider the map 
\begin{eqnarray*}
\prod_\alpha \M_\alpha
& \lto & \M(e,d,(\Lambda_i)_{i=0, \dots, n-1}, (r_\varphi)_\varphi) \\
((\cF^\alpha_i)_i)_\alpha) & \mapsto & \left(\bigoplus_\alpha \cF_i^\alpha\right)_i.
\end{eqnarray*}
associated to these data. All the $\M^\alpha$ are
just isomorphic to $\Spec \cO_E$, so their product is $\Spec \cO_E$ again.

It is not hard to check is that the image of the closed point under
this map is the point described above.

\section{The case of the symplectic group}
\label{symplectic}

Finally, let us consider the question of topological flatness 
of the local model for the symplectic group. 
Let us repeat, with slight notational modifications, the 
definition of the naive local model for the symplectic group given 
in \cite{PR2}; cf. also \cite{RZ}. To simplify the notation, we consider
only the Iwahori case.

Consider a totally ramified extension $F/F_0$ of
degree $e$ as before. Let $V = F^{2g}$ with basis 
 $e_1, \dots, e_g, f_1, \dots, f_g$, and denote by $\{ \cdot, \cdot\}$
the standard symplectic pairing, i.e.
$$ \{e_i, e_j\} = \{f_i, f_j\} = 0, \quad \{e_i, f_{g-j} \} = \delta_{ij}. $$
Let $\delta$ be an $\cO_F$-generator of the inverse different 
 ${\cal D}_{F/F_0}^{-1}$ (if $F$ is tamely ramified over $F_0$, we can take
 $\delta = \pi^{1-e}$).
Let $\langle v, w \rangle = \Tr_{F/F_0}(\delta \{v, w \})$. This
is a non-degenerate alternating form on $V$ with values in $F_0$.

The standard lattice chain $(\Lambda_i)_i$ is self-dual with
respect to $\langle \cdot, \cdot \rangle$.

Now let $r_\varphi = g$ for all $\varphi$.
Associated to $F/F_0$, $V$ and $(r_\varphi)_\varphi$ we have the naive
local model for $GL_{2g}$. 
In this case the reflex field is $F_0$ and $\mu$ is 
just $(e, \dots, e, 0, \dots, 0)$.

The naive local model for the symplectic group is
the closed subscheme of $\M(e, 2g, (\Lambda_i)_i, (r_\varphi)_\varphi)$
consisting of 'self-dual' families of subspaces:
$$\Mn_{GSp} = \{ (\cF_i)_i; \cF_i \lto \Lambda_{i,S} \cong 
\Lambda^\vee_{2g-i, S} \lto \cF^\vee_{2g-i} \text{ is the
zero map} \}. $$
Here $\cdot^\vee$ denotes the $\cO_S$-dual $\hom(\cdot, \cO_S)$,
and the isomorphism $\Lambda_{i,S} \cong \Lambda^\vee_{2g-i, S}$
is the one given by the pairing.

Similarly, we define the local model $\Mloc_{GSp}$ as the
closed subscheme of $\Mloc$ consisting of self-dual families
of subspaces:
$$
\Mloc_{GSp} = \{ (\cF_i)_i \in \Mloc; \cF_i \lto \Lambda_{i,S} \cong 
\Lambda^\vee_{2g-i, S} \lto \cF^\vee_{2g-i} \text{ is the
zero map} \}. $$
Actually, we know that in this case the two local models for
$GL_{2g}$ coincide topologically. Namely, it is enough to check this 
in the maximal parahoric case where the special fibre is irreducible, so 
that one just has to see that the generic and special fibres have the
same dimension. Of course this implies that $\Mn_{GSp}$ and
$\Mloc_{GSp}$ are topologically the same, too. So in order to 
prove the topological flatness, we can work with either one.

\begin{stz}
The local model $\Mloc_{GSp}$ is topologically flat, i.e.
the generic points of the irreducible components of the special fibre
of $\Mloc_{GSp}$ can be lifted to the generic fibre.
The same holds for $\Mn_{GSp}$.
\end{stz}

{\em Proof.}
Since the symplectic local model is defined as a subscheme of the linear
local model, the strata of the special fibre of $\Mloc_{GSp}$ correspond to 
the intersection $\Perm(\mu) \cap \widetilde{W}_{GSp}$ (in the Iwahori case).

The results of Haines and Ng\^o \cite{HN} show that this set coincides with
the $\mu$-admissible set for the symplectic group, and thus the
irreducible components of the special fibre of 
the local model correspond to the
conjugates of $\mu$ under the Weyl group of the symplectic
group (considered as a subgroup of $S_{2g}$). 

In section \ref{parahoric} we generalize Haines' and Ng\^o's results to the
parahoric case. Since this is the only difference between the Iwahori case
and the general parahoric case, until the end of this section we will
assume that we are in the Iwahori case, to simplify notation.

As in the linear case, 
it is enough to show that for each irreducible component there is one 
(suitably chosen) point that can be lifted to the generic fibre.

In this case we cannot hope to find the lifting 
as the image of a map of a product of trivial local models for
the symplectic group, because decomposing the lattice chain as we
did in the linear case would not preserve the pairing.
Nevertheless the proof of the proposition is even simpler
in this case, because the $\mu$ we have to deal with is so special.

We consider the symplectic local model as a closed subscheme
of the linear local model. Let $\nu$ be a conjugate of $\mu$
under the Weyl group of the symplectic group. It is enough to
show that we can lift the point where each $\cF_i$ is given
by the $2ge \times r$-matrix

$$\overline{M} =  \left( \begin{array}{cccc}  
\begin{array}{c} I_{\nu_1} \\
0_{e-\nu_1 \times \nu_1}\end{array} & &  &\\
  &   \begin{array}{c} I_{\nu_2} \\
0_{e-\nu_2 \times \nu_2}\end{array} & & \\
  &  & \ddots & \\
  &  & & \begin{array}{c} I_{\nu_{2g}} \\
0_{e-\nu_{2g} \times \nu_{2g}}\end{array}
\end{array} \right) $$

(which obviously lies inside $\overline{\M}^{sympl}$)
to the generic fibre (of $\M^{sympl}$).

Let us first lift this point to a point in the generic
fibre of the linear local model $\Mloc$. Afterwards we will
show that the lifting actually lies in the symplectic model.

Since all the $\nu_i$ are either $e$ or $0$, we can just lift this 
point by exactly the same matrix over $\cO_E$! It is clear that these
matrices do describe a point in the (linear) local model over $\cO_E$,
in particular that the determinant condition is satisfied.

But since $\nu$ is not an arbitrary permutation of $\mu$, but one under
the Weyl group of the symplectic group, it can never happen that $\nu_i$
and $\nu_{2g-i+1}$ are both 1. In other words, if in the matrix above
there is a unit matrix somewhere in rows $e(i-1)+1, \dots, e(i-1)+e$, 
the columns $e(2g-i) +1, \dots, e(2g-i) + e$ will entirely consist of 0's.
Taking into account that all the $\cF_i$ are given by the same matrix,
it is then clear that this point lies in the symplectic local model.
\qed

We cannot prove that $\M_{GSp_{2g}}$ (or $\Mloc_{GSp_{2g}}$) 
is flat. But since $\M_{GSp_{2g}}$ is topologically flat, we
are led to the

\begin{conj}
The naive local model $\Mn_{GSp_{2g}}$ is flat.
\end{conj}
 
To prove this conjecture, it would be sufficient to show that the 
special fibre is reduced. This would follow if one could identify 
it with a union of Schubert varieties in some affine flag variety,
and by the usual method using Frobenius splittings, it
would be enough to do this for the local models which correspond to 
a maximal parahoric subgroup. But even that seems to be a difficult 
problem in commutative algebra.

\section{The parahoric case}
\label{parahoric}

In this section we will provide
the generalization of the results of Haines and Ng\^o about the
 $\mu$-admissible and the $\mu$-permissible sets
which is needed to prove the topological flatness in the general parahoric
case.

Although we are interested only in the cases of $GL_n$ and $GSp_{2g}$, the
$\mu$-admissible and $\mu$-permissible sets can be defined for any
split connected reductive group. We denote by $\widetilde{W}$ the
extended affine Weyl group, and by $\Omega$ its subgroup of elements of
length zero. In other words, $\Omega$ is the stabilizer of the base alcove.
Then $\widetilde{W}$ is the semi-direct product of the affine Weyl group
 $\Waff$ and $\Omega$.

We fix a dominant coweight $\mu$, and denote by $\tau$
the unique element of $\Omega$ such that $\mu \in \Waff\tau$.

We denote by $\overline{\mathbf a}$ the closure of the base alcove, and by
 $P_\mu$ the convex hull of the translates of $\mu$ under the finite 
Weyl group.

\begin{Def} (Kottwitz-Rapoport)

i) The $\mu$-permissible set $\Perm(\mu)$ is the set of elements $x \in
\Waff\tau$ such that $x(v)-v \in P_\mu$ for all $v \in
\overline{\mathbf a}$.

ii) The $\mu$-admissible set $\Adm(\mu)$ is the set of $x \in
\widetilde{W}$ such that there exists $w \in W_0$ with
 $x \le t_{w\mu}$.

\end{Def}

It was shown by Kottwitz and Rapoport that these two sets coincide for
$G=GL_n$ or $GSp_{2g}$ and minuscule $\mu$, and that the $\mu$-admissible
set is always contained in
the $\mu$-permissible set. Furthermore, we have the following theorem by Haines
and Ng\^o:

\begin{thm} {\rm (\cite{HN}, Theorem 1, Theorem 4, Proposition 5)} 

i) If $\mu$ is a dominant coweight for $GL_n$, then $\Perm(\mu) =
\Adm(\mu)$.

ii) If $\mu$ is a multiple of the dominant minuscule coweight $(1^g, 0^g)$
for $GSp_{2g}$, then $\Perm_{GSp_{2g}}(\mu) = \Perm_{GL_{2g}}(\mu) \cap
 \widetilde{W}_{GSp_{2g}} = \Adm_{GSp_{2g}}(\mu)$.
\end{thm}

(Haines and Ng\^o also show that in general the
 $\mu$-admissible and the $\mu$-permissible set do not coincide.)

It is clear that the set of irreducible components of the special fibre of
$\Mloc$ is exactly the set of maximal elements (with respect to the Bruhat
order) of the $\mu$-permissible set. The theorem says that these maximal
elements are just the conjugates of $\mu$ under the finite Weyl group.

The theorem as it stands relates to the Iwahori case. To prove topological
flatness in the general parahoric case, we need a generalized version which
covers the parahoric case, too. Clearly the Iwahori case is the
most difficult among all the parahoric cases, and as we will show, the
general parahoric case can be deduced from the Iwahori case relatively easily.

\subsection{$GL_n$}

We use the notation of \cite{KR}. Let us recall part of it:

Let $\overline{I}$ be a non-empty subset of $\Z/n\Z$, and denote by $I
\subseteq \Z$ its inverse image under the projection $\Z \lto \Z/n\Z$.
A family $(v_i)_{i\in I}$, $v_i \in \Z^n$, is called a face of type $I$ if
it satisfies the following conditions:
\begin{enumerate}
\item $ v_{i+n} = v_i + {\mathbf 1}$ for all $i \in I$,
\item $v_i \le v_j$ for all $i, j \in I$, $i \le j$,
\item $\sum(v_i) - \sum(v_j) = i-j$ for all $i,j \in I$.
\end{enumerate}

We denote the set of faces of type $I$ by $\faces_I$.

Clearly, a face of type $\Z/n\Z$ is simply an alcove. We denote by $\omega$
the standard alcove $\omega = (\omega_0, \dots, \omega_{n-1})$, $\omega_i =
(1^i, 0^{n-i})$, $i = 0, \dots, n-1$, as well as the corresponding face of
  type $I$, $\omega = (\omega_i)_{i\in I}$.

The extended affine Weyl group $\widetilde{W}$ acts transitively on the set
$\faces_I$. Taking $\omega$ as a base point, we identify $\faces_I$ with
the set $\widetilde{W}/W_I$, where $W_I$ is the stabilizer of 
 $(\omega_i)_{i\in I}$ in $\widetilde{W}$.

If $\overline{J} \subset \overline{I}$ is a non-empty subset, and $J
\subset \Z$ its inverse image under the projection $\Z \lto \Z/n\Z$, we
have a $\widetilde{W}$-equivariant surjection $\faces_I \lto \faces_J$,
defined by $(v_i)_{i\in I} \mapsto (v_i)_{i \in J}$.

We can now adapt the definition above to the parahoric case:

\begin{Def} (Kottwitz-Rapoport)

i) The $\mu$-permissible set $\Perm_I(\mu) \subseteq \faces_I$ is the 
set of elements $(v_i)_{i \in I}$ such that for all $i \in I$,
$v_i - \omega_i \in P_\mu$.

ii) The $\mu$-admissible set $\Adm_I(\mu) \subseteq \faces_I$ is the 
image of $\Adm(\mu)$ under the surjection $\faces_{\Z/n\Z} \lto \faces_I$.
\end{Def}

\begin{stz}
Let $\overline{J} \subseteq \overline{I}$ be a non-empty subset. Then the
restriction of the map $\faces_I \lto \faces_J$
to $\Perm_I$ is a surjection $\Perm_I \lto \Perm_J$.
\end{stz}

{\em Proof.} It is clearly enough to prove the proposition for 
 $\overline{I} = \Z/n\Z$. Given $(v_j)_{j\in J}$ we can then fill in the
 missing $v_j$'s step by step, and thus are reduced to the following lemma.
\qed

\begin{lem} \label{extendPermGLn}
Let $k$, $l$ be integers, $k < l \le k+n$, and let $v_k, v_l \in \Z$.
Assume that $v_k - \omega_k \in P_\mu$, $v_l - \omega_l \in P_\mu$, $v_k -
v_l$ is minuscule and $\sum (v_l) - \sum (v_k) = l - k$.

Then there exists $v_{k+1} \in \Z$ such that $v_{k+1} - \omega_{k+1} \in
P_\mu$, $v_{k+1}-v_k$ and $v_l - v_{k+1}$ are minuscule and 
 $\sum (v_{k+1}) - \sum (v_k) = 1$.
\end{lem}

{\em Proof.} For $v \in \Z^n$ we have $v \in P_\mu$ if and only if 
 $v_{dom} \preceq \mu$, where $v_{dom} \in W.v$ is dominant.

We denote the standard basis vectors of $\Z^n$ by $e_1, \dots, e_n$.
For $m\in \Z$, let $r$ be the unique integer with $1 \le r \le n$,
 $r \equiv m \mod n$. We define $e_m$ to be the basis vector $e_r$.

We are looking for $m$, $1 \le m \le n$, such that $v_{k+1}$ with
$$ v_{k+1} - \omega_{k+1} = (v_k - \omega_k) + e_m - e_{k+1}$$
satisfies the conditions above, i. e. such that
\begin{enumerate}
\item $v_l - v_{k+1}$ is minuscule, i. e. $v_l(m)-v_k(m) = 1$
\item $v_{k+1} - \omega_{k+1} \in P_\mu$
\end{enumerate}
(The other two conditions are satisfied automatically.)

Let $k'$ be the unique integer in $\{1, \dots, n\}$ such that $k' \equiv k+1
(\mod n)$.

{\bf 1st case.} $v_l(k') - v_k(k') = 1$

In this case we can simply choose $m = k'$.

{\bf 2nd case.} $v_l(k') = v_k(k')$

This case is more complicated. Let $\sigma \in S_n$ be such that
 $\sigma(v_k - \omega_k)$ is dominant, and among all such $\sigma$, choose
 $\sigma$ with $\sigma^{-1}(k')$ maximal.

Let
$$ \tilde{\tilde{m}} = \max \{ M; \sigma(v_l)(M) > \sigma(v_k)(M) \}$$
and
$$ \tilde{m} = \min \{ M; \sigma(v_k - \omega_k)(M) 
                           = \sigma(v_k-\omega_k)(\tilde{\tilde{m}}) \}$$

Let $m = \sigma(\tilde{m})$. Obviously this implies $v_l(m)-v_k(m)=1$.
So all that remains to show is that $v_{k+1}-\omega_{k+1} \in P_\mu$ with
this choice of $m$.

By the definition of $\sigma$ and $m$, $\sigma(v_{k+1}-\omega_{k+1})$ is
dominant, too. Thus it is enough to show that $\sigma(v_{k+1}-\omega_{k+1})
\preceq \mu$.

If $\tilde{m} > \sigma^{-1}(k')$, then this is clear. So let us now
consider the case $\tilde{m} < \sigma^{-1}(k')$.

We write $\mu = (\mu(1), \dots, \mu(n)) \in \Z^n$. We have to show that
$$ \sum_{i=1}^N \sigma(v_{k+1}-\omega_{k+1})(i) \le \sum_{i=1}^N \mu(i)$$
for $N = 1, \dots, n-1$. (It is clear that for $N=n$ we have equality since
$\sum_1^n \sigma(v_{k+1}-\omega_{k+1})(i) = 
\sum_1^n \sigma(v_{k}-\omega_{k})(i)$.)

-- For $N \ge \sigma^{-1}(k')$ and for $N < \tilde{m}$, we have
$$ \sum_{i=1}^N \sigma(v_{k+1}-\omega_{k+1})(i) =
   \sum_{i=1}^N \sigma(v_{k}-\omega_{k})(i) \le \sum_{i=1}^N \mu(i),$$
because $v_k - \omega_k \in P_\mu$.

-- Let $\tilde{\tilde{m}} \le N < \sigma^{-1}(k')$.

We have $\sigma(v_k -\omega_k)(k') > \sigma(v_l - \omega_l)(k')$,
and for all $i > \tilde{\tilde{m}}$ we have $\sigma(v_l)(i) =
\sigma(v_k)(i)$ and thus $\sigma(v_l-\omega_l)(i) \le
\sigma(v_k-\omega_k)(i)$. Since furthermore 
 $\sum_{i=1}^n \sigma(v_l-\omega_l)(i) = \sum_{i=1}^n
 \sigma(v_k-\omega_k)(i)$, we see that
$$ \sum_{i=1}^N \sigma(v_k-\omega_k)(i) < \sum_{i=1}^N
\sigma(v_l-\omega_l)(i).$$

So we have
$$ \sum_{i=1}^N \sigma(v_{k+1}-\omega_{k+1})(i) = 
   \sum_{i=1}^N \sigma(v_k-\omega_k)(i) + 1     \le 
   \sum_{i=1}^N \sigma(v_l-\omega_l)(i) \le \sum_{i=1}^N \mu(i).
$$

-- Finally, consider $\tilde{m} \le N < \tilde{\tilde{m}}$.

We know that $\sum_{i=1}^N \sigma(v_k-\omega_k)(i) \le \sum_{i=1}^N
\mu(i)$, since $v_k - \omega_k \in P_\mu$, and we want to show that for
 $\tilde{m} \le N < \tilde{\tilde{m}}$ we even have $<$ here.
For $N = \tilde{\tilde{m}}$ this is certainly true.

Now suppose we had $\sum_{i=1}^N \sigma(v_k-\omega_k)(i) = \sum_{i=1}^N
\mu(i)$ for some $N$, $\tilde{m} \le N < \tilde{\tilde{m}}$.
This implies 
$$ \sum_{i = N+1}^{\tilde{\tilde{m}}} \sigma(v_k - \omega_k)(i) <
\sum_{i=\N+1}^{\tilde{\tilde{m}}} \mu(i),$$
and thus $\sigma(v_k - \omega_k)(\tilde{m}) < \mu(N+1) \le \mu(N)$ (because
for $\tilde{m} \le N < \tilde{\tilde{m}}$ all $\sigma(v_k-\omega_k)(N)$ are
equal).  But then we get
$$ \sum_{i=1}^{N-1} \sigma(v_k-\omega_k)(i) > \sum_{i=1}^{N-1} \mu(i),$$
which is a contradiction.

The lemma is proved.\qed

\begin{kor}
Let $\mu$ be a dominant coweight for $GL_n$, and let $I$ be the inverse
image of a non-empty subset $\overline{I} \subseteq \Z/n\Z$.
Then $\Perm_I(\mu) = \Adm_I(\mu)$.
\end{kor}

{\em Proof. }
The $\mu$-admissible set is always contained in the $\mu$-permissible
set. Since we know that in the Iwahori case the two sets coincide, and
because we have surjections $\Perm(\mu) \lto \Perm_I(\mu)$ (by the
proposition) and $\Adm(\mu) \lto \Adm_I(\mu)$ (obvious), the corollary
follows.
\qed

\subsection{$GSp_{2g}$}

Now let $G = GSp_{2g}$. Since the proofs for the symplectic group are based
on reductions to the linear case, we use a subscript $\cdot_G$ to denote
data corresponding to the symplectic group; notation without subscript
refers to the $GL_{2g}$-case, as in the previous section.

We denote by $\Theta: \Z^{2g} \lto \Z^{2g}$ the automorphism given by
$(x_1, \dots, x_{2g}) \mapsto (-x_{2g}, \dots, -x_1)$.

This automorphism acts on the root system of $GL_{2g}$, and the
'$\Theta$-invariant part', denoted $R^{[\Theta]}$, is the root system
of $GSp_{2g}$. See \cite{HN}, \S\S 9, 10. In particular, the extended affine
Weyl group $\widetilde{W}_G$ for the general symplectic group is a subgroup
of the extended affine Weyl group for the general linear
group. Furthermore, by \cite{HN} Proposition 9.6, the Bruhat order on 
 $\widetilde{W}_G$ is inherited from the Bruhat order on $\widetilde{W}$.

The vectors $\eta_i = \frac{1}{2} (\omega_i + \omega_{2g-i})$, 
 $i = 0, \dots, g$,
serve as 'vertices' of the base alcove for the symplectic group. 

Let $\mu$ be a dominant coweight for $G$.
We can consider the $\mu$-admissible set $\Adm_G(\mu)$ and the
 $\mu$-permissible set $\Perm_G(\mu)$ as
subsets of $\widetilde{W}$. We have
\begin{eqnarray*}
\Adm_G(\mu) & = & \{ x \in \widetilde{W}; \quad
                     x \le t_{w\mu} \text{ for some } w \in W_{0,G}\}, \\
\Perm_G(\mu) & = & \{ x \in \widetilde{W}; \quad
                      x(\eta_i) - \eta_i \in P_{G, \mu} \text{ for all } 
                      i = 0, \dots, g \} .
\end{eqnarray*}

Here we denote by $P_{G,\mu}$ the convex hull of the translates of $\mu$ under the
finite Weyl group (of $Sp_{2g}$).

Now let $\overline{I} \subseteq \Z/2g\Z$ be a non-empty symmetric subset,
i.e. a non-empty subset such that
its inverse image $I$ under the projection $\Z \lto \Z/2g\Z$ satisfies
 $I = -I$. (Obviously there is a one-to-one correspondence between
 symmetric subsets of $\Z/2g\Z$ and subsets of $\{0, \dots, g \}$.)  
We describe the set of $G$-faces of type $I$ as a subset of $\faces_I$:
A $G$-face of type $I$ is a face $(v_i)_{i\in I}$ for
$GL_{2g}$ such that there exists $d \in \Z$ with
$$ v_{2g-i} = \Theta(v_i) + (d^{2g})$$
for all $i \in I$. Clearly a $G$-face of type $\Z/2g\Z$ is simply an alcove.
We denote the set of $G$-faces of type $I$ by $\faces_{G,I}$.

We obtain a commutative diagram
\begin{diagram}
\widetilde{W}_G / W_{G,I} & \rTo^\cong & \faces_{G, I} \\
\dInto & & \dInto \\
\widetilde{W} / W_{I} & \rTo^\cong & \faces_I.
\end{diagram}

Here $W_{G,I} =  W_I \cap \widetilde{W}_G \subseteq \widetilde{W}_G$ 
is the stabilizer of $(\omega_i)_{i \in I}$ in $\widetilde{W}_G$.

We can now define parahoric versions of the admissible and permissible
sets.

\begin{Def}
i) The $\mu$-permissible set $\Perm_{G,I}(\mu) \subseteq 
\widetilde{W}_G / \widetilde{W}_{G,I} \cong \faces_{G,I}$ is the 
set of elements $x\widetilde{W}_{G,I}$ such that for all 
 $i \in I \cap \{0, \dots, g\}$, $x(\eta_i) - \eta_i \in P_{G,\mu}$.

ii) The $\mu$-admissible set $\Adm_{G, I}(\mu) \subseteq \faces_{G,I}$ is the 
image of $\Adm_G(\mu)$ under the surjection $\faces_{G, \Z/2g\Z} \lto
\faces_{G, I}$.
\end{Def}

It turns out, however, that the set which naturally describes the strata
occuring in the special fibre of a local model is not the $\mu$-permissible
set but the intersection $\Perm_I(\mu) \cap \widetilde{W}_G$. The goal of
this section is to show that actually 
$$ \Adm_{G, I}(\mu) = \Perm_{G, I}(\mu) = \Perm_I(\mu) \cap
\widetilde{W}_G. $$

The key point is to show that the natural map
 $ \Perm(\mu) \cap \widetilde{W}_G \lto Perm_I(\mu) \cap
 \widetilde{W}_G/\widetilde{W}_{G,I}$ is surjective.

If $\pi \colon \bar{J} \subseteq \bar{I}$ is a non-empty symmetric subset,
we have a
$\widetilde{W}_G$-equivariant surjection $\faces_{G,I} \lto \faces_{G,J}$.

Now suppose $\bar{J} \subset \bar{I} \subseteq \Z/2g\Z$ are non-empty
symmetric subsets, such that $\bar{I} = \bar{J} \cup \{ \bar{k} + \bar{1},
 -(\bar{k} + \bar{1}) \}$ for some $k\in J$ with $k+1 \not\in J$. (For an
 integer $m$ we denote by $\bar{m}$ its class in $\Z/2g\Z$.)
Let $l$ be the smallest integer in $J$ which is greater than $k$.

\begin{lem} {\rm (\cite{KR}, Lemma 10.3)}
In the situation above we have a bijection
$$\pi^{-1}({\mathbf v}) \lto \left\{ w \in \Z^{2g}; v_k \le w \le v_l,
\sum(w) = \sum(v_k) + 1 \right\}, (w_i)_{i\in I} \mapsto w_{k+1}.$$
\end{lem}

\begin{stz}
Let $\mu$ be a positive multiple of the dominant minuscule coweight
 $(1^g, 0^g)$ for
$G$.
Let $\bar{J} \subseteq \bar{I}$ be a non-empty symmetric subset. Then the
natural map $\Perm_I(\mu) \cap \widetilde{W}_G/W_{G,I} \lto \Perm_J(\mu) \cap
\widetilde{W}_G/W_{G,I}$ is surjective.
\end{stz}

{\em Proof.}
It is clearly enough to consider $\bar{J} \subseteq \bar{I}$ as in the
lemma. Let $(v_j)_{j\in J} \in \Perm(\mu)_{G,J}$. We would like to extend this
face to a $G$-face of type $I$ by defining suitable $v_{k+1}, v_{-(k+1)}
\in \Z^{2g}$. By the lemma above and lemma \ref{extendPermGLn}, we find a
$G$-face $(w_i)_{i\in I}$ of type $I$ which maps to $(v_j)_j$ under $\pi$
and such that $w_{k+1} - \omega_{k+1} \in P_{G,\mu}$.

We have to show that then $w_{-(k+1)} - \omega_{-(k+1)}$ automatically
holds, too.
Now $\mu$ is of the form $(d^g, 0^g)$ for some $d$, and 
$$ w_{-(k+1)} - \omega_{-(k+1)} = (d^{2g}) + \Theta(w_{k+1} - \omega_{k+1}) $$
since $(w_i)_{i\in I}$ is a $G$-face. 

Since $\mu$ has this special form, for a dominant coweight $\lambda$ 
we have $\lambda_{dom} \preceq \mu$ if and only if $\lambda(i) \le d$ for
all $i$ (and $\sum(\lambda) = gd$). If this holds for some $\lambda$,
it clearly holds for $(d^{2g}) -\lambda$ as well, so we are done.
\qed

\begin{kor}
Let $\mu$ be a positive multiple of the dominant coweight for $G$, 
and let $\overline{I} \subseteq \Z/2g\Z$ be a non-empty symmetric subset.
Then $\Perm_{G,I}(\mu) = \Adm_{G,I}(\mu)$.
\end{kor}

{\em Proof.}
First, we have an inclusion $\Perm_{G, I}(\mu) = \Perm_I(\mu) \cap
\widetilde{W}_G/W_{G,I}$. This is theorem 10.1 in \cite{HN} for $I=\Z/2g\Z$,
and it is easy to see that this is a 'vertex-by-vertex' proof, i.e. it
works for arbitrary I.

Now recall that in addition we know that in any case the $\mu$-admissible set
is contained in the $\mu$-permissible set.

Since $\Adm_{G,I}(\mu)$ and $\Perm_I(\mu) \cap
\widetilde{W}_G/W_{G,I}$ both coincide with the image of 
 $\Adm_G(\mu) = \Perm(\mu) \cap \widetilde{W}_G$, the corollary follows.
\qed

\end{document}